\documentclass[english,a4paper,11pt]{article}
\usepackage{theorem}
\usepackage{amssymb,amsmath}
\usepackage{latexsym}
\usepackage{babel}
\usepackage{color}
\usepackage{slashbox}

\usepackage{color,graphicx,array,fleqn}
\usepackage{bbm}

\newtheorem{theorem}{Theorem}[section]

\newtheorem{corollary}[theorem]{Corollary}

\theorembodyfont{\normalfont\mdseries\upshape}
\newtheorem{remark}[theorem]{Remark}

\newtheorem{example}[theorem]{Example}

\title{\bf A minimax approach for inverse variational inequalities}
\author{P. Montiel L\'opez}

\begin{document}

\maketitle

\centerline{University of Granada, Centro de Magisterio La Inmaculada,} 

\centerline{Department of Sciences, Granada (Spain), e-mail: pablomontiel@eulainmaculada.com}

\begin{abstract}
\noindent In this work, we characterize the existence of a solution for a certain variational inequality by means of a classical minimax theorem. In addition, we propose a numerical algorithm for the solution of an inverse problem associated with a variational inequality. To this end we state a collage-type result in this variational framework. 
\end{abstract}



\section{Introduction}

The well-known \textsl{collage theorem} \cite{barervharlan86}, a direct consequence of the Banach fixed point theorem, states that for the unique fixed point $x_0$ of a $c$-contractive self-mapping $\Psi$ on a complete metric space $(X,d)$, there holds
\[
x \in X \ \Rightarrow \ d(x,x_0) \le \frac{1}{1-c}d(x,\Psi(x)).
\]
This result has became a fundamental tool for a class of numerical methods providing the solution of some inverse problems: see, for instance, \cite{kungom03,kunmenlatvrs12,kunvrs99}. Furthermore, if one replaces the Banach fixed point theorem by the Lax--Milgram theorem or some of its generalizations (\cite{garrui14,rui09}), it is possible to establish results along the lines of the collage theorem: these are the so-called \textsl{generalized collage theorems} (\cite{berkunlatrui16,kunlatlevrui15,kunlatvrs09,kunmenlatvrs12}). Our main aim in this paper is to state a collage-type result starting from the Stampacchia theorem to deal with inverse problems related to a variational inequality. 

First of all, we prove that the existence of a solution for a variational inequality in a reflexive Banach space is equivalent to the existence of a constant satisfying an adequate convexity condition. To this end, we make use of the classical minimax theorem of J. von Neumann and K. Fan. In addition, we show how our result clearly implies the Stampacchia theorem. Once the existence of a solution is studied, we deal with the inverse problem associated with a variational inequality. The above-mentioned collage theorem --a stability result derived from the Stampacchia theorem-- and the use of an adequate Schauder basis in the involved reflexive Banach space, allow us to design a numerical method for the solution of the variational inverse problem. We illustrate our results with a numerical test.

Let us recall some standard notations. For a real normed space $E$, $B_E$ denotes its \textsl{closed unit ball} and $E^*$ its \textsl{topological dual}, that is, the Banach space of those continuous and linear functionals on $E$. If $m \ge 1$, $\Delta_m$ stands for the \textsl{unit simplex} of $\mathbb{R}^m$, that is,
\[
\Delta_m:= \{t=(t_1,\dots,t_m) \in \mathbb{R}^m : \ t_1,\dots,t_m \ge 0, \ t_1+\cdots +t_m=1\}.
\]

\section{Existence of a solution for variational inequalities}

We adopt a minimax approach for deaing with variational inequalities, as in \cite{garrui19,rui14}, unlike that of K. Fan \cite{fan72} (see also \cite{aub98}), where an equilibrium result is the main tool. In this way, we not only derive sufficient but also necessary conditions for the existence of a solution for a certain variational inequality.

We make use of the following classical minimax theorem (\cite{fan53,kaskol96,neu28}). More general or different versions can be found, for instance, in \cite{pat17,ric17,rui18,sai18,simrui02,syg18,tia17}.

\bigskip

\begin{theorem}\label{th:minimax}
Assume that $X$ is a nonempty, convex and compact subset of a real topological vector space, $Y$ is a nonempty set and $g: X \times Y \longrightarrow \mathbb{R}$ is continuous and concave on $X$. Then, 
\[
\hbox{there exists } x_0 \in X : \ \inf_{y \in Y} \max_{x \in X} g(x,y) \le \inf_{y \in Y} g(x_0,y)
\]   
if (and only if)
\[
\left.
	\begin{array}{c}
	m \ge 1, \ t \in \Delta_m  \\
	y_1, \dots,y_m \in Y
\end{array}	  
\right\}
\ \Rightarrow \
\inf_{y \in Y} \max_{x \in X} g(x,y) \le \max_{x \in X} \sum_{j=1}^m t_j g(x,y_j).
\]
\end{theorem}

\bigskip

This result has been used, in an equivalent form of theorem of the alternative, or that of Hahn--Banach type result, to characterize the existence of a solution for nonlinear infinite programs: see  \cite{monrui19,monrui17,monrui16,rui16-1}. Let us also mention that the convexity condition in Theorem \ref{th:minimax} is the so-called \textsl{infsup}-convexity (see \cite{kaskol96,rui16-2}).

Now we state a general theorem of existence for certain variational inequalities that implies the classical Stampacchia theorem (\cite{sta64}).

\bigskip

\begin{theorem}\label{th:existence}
Let $E$ be a real reflexive Banach space, $x_0^* \in E^*$, $a: E \times E \longrightarrow \mathbb{R}$ be a continuous bilinear form, and $Y$ be a nonempty weak closed subset of $E$. Then
\begin{equation}\label{eq:47}
\hbox{there exists } x_0 \in Y : \ y \in Y \ \Rightarrow \ x_0^*(y-x_0) \le a(y,y-x_0)
\end{equation}
if, and only if, for some $\alpha \ge 0$, $Y \cap \alpha B_E \neq \emptyset$ and
\begin{equation}\label{eq:47bis}
\left.
	\begin{array}{c}
	m \ge 1, \ t \in \Delta_m  \\
	y_1, \dots,y_m \in Y
\end{array}	  
\right\}
\ \Rightarrow \
\sum_{j=1}^m t_j (x_0^*(y_j)-a(y_j,y_j))  \le \max_{x \in Y \cap \alpha {B_E}} \left(x_0^*(x)-a\left(\sum_{j=1}^m t_j y_j,x\right) \right).
\end{equation}
\end{theorem}

\noindent \textsc{Proof}. The fact that \eqref{eq:47} implies \eqref{eq:47bis} is very easy to check: it suffices to consider $\alpha:=\|x_0\|$ and use the linearity of $a$ at its first variable.

And conversely. Let $X:=Y \cap \alpha B_E$ and observe that, thanks to \eqref{eq:47bis} with $m=1$, it holds that
\[
\begin{array}{rl}
0 & \le \displaystyle \inf_{y \in Y} \left( a(y,y)-x_0^*(y)- \max_{x \in X} (x_0^*(x)-a(y,x))\right)   \\
  & = \displaystyle \inf_{y \in Y} \max_{x \in X} (a(y,y-x)-x_0^*(y-x)).
\end{array}
\]
Then, we write
\[
\mu :=  \inf_{y \in Y} \max_{x \in X} (a(y,y-x)-x_0^*(y-x))
\]
and so, $\mu \in \mathbb{R}$. If we define the bifunction $f: X \times Y \longrightarrow \mathbb{R}$ by
\[
f(x,y):=a(y,y-x)-x_0^*(y-x)-\mu, \qquad (x \in X, \ y \in Y),
\]
then
\[
0 =\inf_{y \in Y} \max_{x \in X} f(x,y),
\]
and the minimax theorem, Theorem \ref{th:minimax}, applies when we endow $E$ with its weak topology, since $X$ is a nonempty, weak compact and convex subset of $E$ ($Y$ is weak closed and $B_E$ is weak compact) and $f$ is concave and weak continuous on $X$. Therefore, there exists $x_0 \in X$ such that
\[
\inf_{y \in Y} \max_{x \in X} f(x,y) \le \inf_{y \in Y} f(x_0,y),
\]
i.e., there holds \eqref{eq:47}, if
\[
\left.
	\begin{array}{c}
	m \ge 1, \ t \in \Delta_m  \\
	y_1, \dots,y_m \in Y
\end{array}	  
\right\}
\ \Rightarrow \
0 \le \max_{x \in X} \sum_{j=1}^m t_j f(x,y_j),
\]
which is clearly equivalent to the condition \eqref{eq:47bis}.
\hfill${\Box}$

\bigskip

When, in addition, $Y$ is convex --in this case $Y$ is closed for the norm topology, since it is weakly closed and convex-- and $a$ is non-negative at the diagonal, the condition \eqref{eq:47bis} is simpler and \eqref{eq:47} can be equivalently reformulated:

\bigskip

\begin{corollary}
Under the same assumptions that in Theorem \ref{th:existence}, let us also suppose that $Y$ is convex and that, for all $x \in E$, $a(x,x) \ge 0$. Then
\begin{equation}\label{eq:46}
\hbox{there exists } x_0 \in Y : \ y \in Y \ \Rightarrow \ x_0^*(y-x_0) \le a(x_0,y-x_0)
\end{equation}
if, and only if, there exists $\alpha \ge 0$ such that $Y \cap \alpha B_E \neq \emptyset$ and
\begin{equation}\label{eq:46bis}
y \in Y \ \Rightarrow \ x_0^*(y)-a(y,y) \le \sup_{x \in Y \cap \alpha{B_E}} (x_0^*(x)-a(y,x)).
\end{equation}
\end{corollary}

\noindent \textsc{Proof} According to \cite[Lemma 4.1]{cap14} and the positivity condition on the diagonal of $a$, the variational problem \eqref{eq:46} is equivalent to \eqref{eq:47}. Therefore, in view of Theorem \ref{th:existence}, we will state the equivalence of \eqref{eq:46bis} and \eqref{eq:47bis}. 

It is clear that \eqref{eq:47bis} $\Rightarrow$ \eqref{eq:46bis}, hence we focus on proving the converse. So, let $m \ge 1$, $t \in \Delta_m$ and $y_1,\dots,y_m \in Y$. Then, the convexity of $Y$, the linearity of $f$ and the convexity of the quadratic form $x \in E \mapsto a(x,x)$ (taking into account that $a$ is non-negative on its diagonal) yield
\[
\begin{array}{rl}
	\displaystyle \sum_{j=1}^m t_j (x_0^*(y_j)-a(y_j,y_j)) & = \displaystyle  \sum_{j=1}^m t_j x_0^*(y_j) - \sum_{j=1}^m t_j a(y_j,y_j)   \\
	& \le \displaystyle x_0^* \left( \sum_{j=1}^m t_j y_j \right) -a \left(  \sum_{j=1}^m t_j y_j ,  \sum_{j=1}^m t_j y_j \right)   \\
	& \displaystyle \le \sup_{x \in Y \cap \alpha {B_E}} \left( x_0^*(x)- a \left(  \sum_{j=1}^m t_j y_j , x \right) \right).
\end{array}
\]
\hfill${\Box}$

\bigskip

The existence of solution for (systems of) variational equations or different variational inequalities has been previously established by means of minimax inequalities or the Hahn--Banach theorem \cite{garrui19,garrui14,rui18,rui14,rui09}.

Let us note that Stampacchia's theorem is a direct consequence of the previous result: assume that $E$ is a Hilbert space, $Y$ is a nonempty, closed and convex subset of $E$,  $a: E \times E \longrightarrow \mathbb{R}$ is $\rho$-coercive ($\rho > 0$), continuous and bilinear  and that $x_0^* \in E^*$ --here we do not use the Riesz identification of $E$ with its topololgical dual $E^*$--. These hypotheses are sufficient for the existence of $x_0 \in E$ with
\[
y \in Y \ \Rightarrow \ x_0^*(y-x_0) \le a(x_0,y-x_0).
\]
To prove this, let $\beta >0$ such that $Y \cap \beta B_E \neq \emptyset$ and note that
\[
\begin{array}{rl}
\displaystyle \frac{x_0^*(y)-a(y,y)}{\| y \|}-\frac{\displaystyle \sup_{x \in Y \cap \beta {B_E}}(x_0^*(x)-a(y,x))}{\| y \|} & \le \|x_0^*\|-\rho \| y \|+ \displaystyle \beta\frac{\|x_0^*-a(y,\cdot)\|}{\|y\|}   \\
 & \le \| x_0^* \| \left( 1 + \displaystyle \frac{\beta}{\| y \|}  \right)+ \beta \| a \|-\rho \| y \|, 
\end{array}
\]
hence in particular, for some $\alpha > \beta$ we have that ($Y \cap \alpha B_E \neq \emptyset$ and)
\[
y \in Y, \ \| y \| > \alpha \ \Rightarrow \ x_0^*(y)-a(y,y) \le \sup_{x \in Y \cap \alpha {B_E}} (x_0^*(x)-a(y,x)),
\]
while, trivially,
\[
y \in Y \cap \alpha B_E \ \Rightarrow \ x_0^*(y)-a(y,y) \le \sup_{x \in Y \cap {B_E}} (x_0^*(x)-a(y,x)),
\]
so condition \eqref{eq:46bis} is valid and then the variational problem \eqref{eq:46} admits a solution.

\section{The inverse problem}

The Stampacchia theorem straightforwardly implies the following collage-type result, which will allow us to deal with an inverse problem associated with a certain variational inequality. The role of this result is the same as that of the Banach fixed point in the collage treatment of some inverse problems (see, for instance, \cite{kungom03}). In fact, such an approach is motivated by this stability idea, previously considered for differential and variational equations in terms of the Lax--Milgram theorem or some of its generalizations, which is known as the \textsl{generalized collage theorem}: see \cite{berkunlatrui16,kunlatlevrui15,kunlatvrs09,kunmenlatvrs12,kunvrs99}.

\bigskip

\begin{theorem}\label{th:collage}
Let $J$ be a nonempty set, $E$ be a real Hilbert space, $Y$ be a nonempty, closed and convex subset of $E$, and for all $j \in J$ let $x_j^* \in E^*$, $a_j: E \times E \longrightarrow \mathbb{R}$ be continuous, bilinear and such that, for some $\rho_j >0$,
\[
y \in E \ \Rightarrow \ \rho_j \| y \|^2 \le a_j(y,y). 
\]
If in addition $x_j\in Y$ is the solution of the variational inequality 
\[
y \in Y \ \Rightarrow \ x_j^*(y-x_j) \le a(x_j,y-x_j),
\]
then, 
\[
y \in Y, \ j \in J \  \Rightarrow \ \|y-x_j \| \le \frac{\|a_j(y,\cdot)-x_j^*\|}{\rho_j}.
\]
\end{theorem}

\noindent \textsc{Proof} Given $y \in Y$ and $j \in J$, the announced inequality follows from this chain of inequalities:
\[
\begin{array}{rl}
\rho \| y-x_j \|^2 & \le a_j(y-x_j,y-x_j)   \\
                   & = a_j(y,y-x_j)-a_j(x_j,y-x_j)   \\
                   & \le a_j(y,y-x_j)-x_j^*(y-x_j)   \\
                   & = (a_j(y,\cdot)-x_j^*)(y-x_j)   \\
                   & \le \| a_j(y,\cdot)-x_j^*\| \|y-x_j\|.
\end{array}
\]
\hfill${\Box}$

\bigskip

\begin{remark}\label{re:norm} We point out that, if $Y$ is a closed affine subspace of $E$, then we can replace in the theorem $\|a_j(y,\cdot)-x_j^*\|$ by
$\|(a_j(y,\cdot)-x_j^*)_{| Y_0}\|$, where $Y_0$ is the closed linear subspace $Y-Y$ of $E$.

\end{remark}

\bigskip

As mentioned above, our work is motivated by the collage treatment in \cite{kunlatvrs09}, here for the following inverse problem: let us assume the hypotheses in Theorem \ref{th:collage} hold. If $y\in Y$ is a given target element, then we want to determine, if possible, that $j_0 \in J$ minimizing the distance $\| y-x_j\|$, that is, 
\[
\| y-x_{j_0}\|=\inf_{j \in J} \| y-x_j\|.
\]
However, such an optimization problem is very difficult to solve, since we  should not only minimize that function but also solve all the variational problems: find $x_j \in Y$ with 
\begin{equation}\label{eq:inequality}
y \in Y \ \Rightarrow \ x_j^*(y-x_j) \le a(x_j,y-x_j).
\end{equation}
If, in addition, the family of continuous and bilinear forms is uniformly coercive, in the sense that
\[
\rho:= \inf_{j \in J} \rho_j >0,
\]
then, in view of this assumption, Theorem \ref{th:collage} and Remark \ref{re:norm}, if $Y$ is a closed affine subspace of $E$ and $y \in Y$ is a given target element, we can replace that optimization problem by this other one:
\begin{equation}\label{eq:minimun}
\inf_{j \in J} \|(a_j(y,\cdot)-x_j^*)_{| Y_0}\|.
\end{equation}
Let us note that the data of this nonlinear program only depend on the data of the variational inequalities \eqref{eq:inequality}. Moreover, this new optimization problem 
can be approximately solved by means of an orthogonal basis in the Hilbert space $E$ (\cite{kunlatvrs09}), or even by a Schauder basis in $E$ (\cite{berkunlatrui16,kunlatlevrui15,kunmenlatvrs12,kunvrs99}). Let us notice that Schauder bases are tools for the numerical solution of a wide variety of
differential, integral and integro-differential problems \cite{berfergarrui09,berforgarrui04, gamgarrui05,gamgarrui09,palrui05}.
It is worth mentioning that, unlike in our motivating works \cite{berkunlatrui16,kunlatlevrui15,kunlatvrs09}, the choice of the target function is obtained by an approximation of the exact solution provided by a Galerkin scheme for the direct problem. Therefore, such a numerical method  will become an auxiliary tool for dealing with the inverse problem.

Now, we illustrate this collage-based numerical method  for the solution of the inverse problem associated with the family of variational inequalities \eqref{eq:inequality}. Before this, we introduce the above-mentioned numerical method for approximating the direct problem.

\bigskip

\begin{example}\label{ex:1}
Consider the boundary value problem
\[
\left\{\begin{array}{ll}

-u''(x)+j u(x)=f(x) \quad \hbox{on } (0,1)\\

u(0)=\alpha, \ u(1)=\beta,
\end{array}
\right. 
\]
where $\alpha, \beta , j \in \mathbb{R}$ with $j >0$, and $f \in L^\infty (0,1)$. Standard reasoning leads to the following variational formulation of this nonhomgeneous problem: let $v$ be a test function in the closed and convex subset $Y$ of $H^1 (0,1)$
\[
Y:=\left\lbrace v \in H^1(0,1): \ v(0)=\alpha , \ v(1)=\beta  \right\rbrace,
\]
multiply the second order differential equation by $v-u$ and integrate by parts to arrive at
\begin{equation}\label{eq:1}
v\in Y \ \Rightarrow \ \int^1_0 u'(v - u)'+ j \int^1_0 u(v - u) \geq \int^1_0 f(v - u).
\end{equation}
Then, the Stampacchia theorem, when applied to the continuous, bilinear and $1$-coercitive form $a_j:H^1(0,1) \times H^1(0,1) \longrightarrow \mathbb{R}$ 
\[
a_j(u,v):=\int^1_0 u'v'+ j \int^1_0 uv, \qquad (u,v \in H^1(0,1))
\]
and the continuous and lineal functional $x_0^*: H^1(0,1) \longrightarrow \mathbb{R}$
\[
x_0^*(v):= \int^1_0 fv, \qquad (v \in H^1(0,1))
\]
guarantees the existence of a unique solution $u \in Y$ of the variational inequality \eqref{eq:1}.

Now we introduce a Galerkin method based upon the properties of a certain Schauder basis in $H^1_0(0,1)$. To this end, if we write \eqref{eq:1} with $\omega \in H^1_0 (0,1)$,  $\omega:=v-u$ we obtain
\[
\omega  \in H^1_0 (0,1) \ \Rightarrow \  \int^1_0 u'w'+ j \int^1_0 uw= \int^1_0 fw.
\]
In order to generate an increasing sequence of finite dimensional linear subspaces of $H^1_0(0,1)$ whose union is dense in this space, let us consider the Haar $\left\lbrace h_k \right\rbrace_{k \geq 1}$ in $L^2 (0,1)$ and define the sequence in $H^1(0,1)$
\[
g_1(x):=1, \qquad (x \in [0,1])
\]
and for any $k \ge 2$,
\[
g_k(x):=\int^x_0 h_{k-1} (t) dt, \qquad (x \in [0,1]).
\]
This sequence is a Schauder basis for $H^1(0,1)$ and the sequence $\{g_{k+2}\}_{k \ge 1}$ is also a Schauder basis for $H^1_0(0,1)$ (see 
\cite[Propositions 4.7 and 4.8]{fu}. Let us note that we can express any element $v \in Y$ as
\[
v= \alpha + (\beta - \alpha)x+  \sum^{\infty}_{k=1} \alpha_{k} g_{k+2},
\]
for some scalars $\alpha_k$. Then, we define the aforementioned $m$-dimensional subspaces
\[
H_m:=\mathrm{span}\left\lbrace g_3, \dots, g_{m+2}  \right\rbrace, \qquad  (m \ge 1).
\]
In the following table we collect the errors generated when approximating the exact solution $u$ of the variational inequality \eqref{eq:1} by the solution $u_m \in Y$ of the $m$-dimensional problem, with the data $\alpha=-3$, $\beta=-4$, $j =\sqrt{2}$,
\[
f(x):= -2+\sqrt{2}(x^2 -2x -3), \qquad (x \in [0,1])
\]    
and $m=3,7,15,31,63$.

\begin{table}[h]
\begin{center}
\begin{tabular}{|l|l|l|l|}
\hline
$m$ & $\Vert u - u_{m} \Vert_{L_{2}}  $&   $\Vert u' - u'_m \Vert_{L_{2}} $ &  $ \Vert u - u_m \Vert_{H_{1}}  $ \\
\hline 
$3$ & $0.0105048 $ & $0.144383 $ & $ 0.144765$ \\ \hline
$7$ & $0.00261572 $ & $0.0721747 $ & $ 0.0722221 $ \\ \hline
$15$ & $0.000653279 $ & $ 0.0360851$  & $ 0.0360911$\\ \hline
$31$ & $0.000163279 $ & $0.0180423 $  & $0.018043 $\\ \hline
$63$ & $0.0000408172 $ & $0.00902111 $  & $0.0090212 $\\ \hline
\end{tabular}
\caption{Errors of the direct method.}
\label{tabla:sencilla}
\end{center}
\end{table}

\end{example}

\bigskip

Now we can address the inverse problem:

\bigskip

\begin{example}
Let us now introduce the following inverse problem associated with  the variational equation \eqref{eq:1}: 
\[
\left\{\begin{array}{ll}
-u''(x)+ j u(x)=f(x) \quad \hbox{on } (0,1) \\
u(0)=-3, \ u(1)=-4,
\end{array}
\right. 
\]
where $1\leq j \leq 4$ and 
\[
f(x):= -2+\sqrt{2}(x^2 -2x -3), \qquad (x \in [0,1]).
\]    
We analyze the performance of the collage-based method. In order to illustrate it, we take $j =\sqrt{2}$ and obtain a target function $u_m \in Y$ by the Galerkin method described in Example \ref{ex:1}. Then, in order to solve the inverse problem given in  (\ref{eq:minimun}), we observe that in our problem
\[
\inf_{j \in [1,4]} \|a_{j}(y,\cdot)-x^*\|=\inf_{j \in [1,4]} \sup_{\scriptsize{\begin{array}{cc}\omega \in H^1_0(0,1)\\
\|\omega\|=1\end{array}}} |a_{j}(y,\omega)-x^*(\omega)|.
\]
Now, we note that $\omega \in H^1_0(0,1)$ can be written as
\[
\omega=\sum_{k=1}^\infty \alpha_k g_{k+2},
\]
for some scalars $\alpha_k$ verifying that $|\alpha_k|\leq 1$. Then
\[
\inf_{j \in [1,4]}  \sup_{\scriptsize{\begin{array}{cc}\omega \in H^1_0(0,1)\\
\|\omega\|=1\end{array}}} |a_{j}(y,\omega)-x^*(\omega)| \leq \inf_{j \in [1,4]}  \left|\sum_{k=1}^{\infty}(a_{j}(y,g_{k+2})-x^*(g_{k+2})\right|.
\]

The next table shows the approximations obtained for $j$ solving the minimizing problem
\[
 \inf_{j \in [1,4]}  \left|\sum_{k=1}^{n}(a_{j}(y_m,g_{k+2})-x^*(g_{k+2})\right|.
\]
for the value $n=31,$ with different target elements $u_m$ of the previous example.

\begin{table}[h]
\begin{center}
\begin{tabular}{|c|l|l|l|l|}
\hline
\backslashbox{$n$}{$m$} & $3$ &   $ 7 $ &  $ 15  $ & $31$  \\
\hline  \hline
$31$ & $1.53389 $ & $1.46679 $ & $ 1.43170$ & $1.41421$  \\ \hline \hline

\end{tabular}
\caption{Numerical results for the inverse problem.}
\label{tabla:sencilla}
\end{center}
\end{table}

\end{example}

\end{document}